\documentclass[12pt,english,reqno]{amsart}
\usepackage{amssymb,amsmath,amstext,amsthm,amsfonts}
\usepackage[dvips]{graphicx}
\usepackage[ansinew]{inputenc}

\newtheorem{theorem}{Theorem}

\newtheorem{maintheorem}{Theorem}

\newcommand{\cmt}{\begin{maintheorem}}
\newcommand{\fmt}{\end{maintheorem}}

\newtheorem{maincorollary}[maintheorem]{Corollary}

\newcommand{\cmc}{\begin{maincorollary}}
\newcommand{\fmc}{\end{maincorollary}}

\newtheorem{T}{Theorem}[section]
\newcommand{\cte}{\begin{T}}
\newcommand{\fte}{\end{T}}

\newtheorem{Corollary}[T]{Corollary}
\newcommand{\cco}{\begin{Corollary}}
\newcommand{\fco}{\end{Corollary}}

\newtheorem{Proposition}[T]{Proposition}
\newcommand{\cpr}{\begin{Proposition}}
\newcommand{\fpr}{\end{Proposition}}

\newtheorem{Lemma}[T]{Lemma}
\newcommand{\cle}{\begin{Lemma}}
\newcommand{\fle}{\end{Lemma}}

\newcommand{\csle}{\begin{Lemma}}
\newcommand{\fsle}{\end{Lemma}}

\newtheorem{Remark}[T]{Remark}
\newcommand{\cre}{\begin{Remark}}
\newcommand{\fre}{\end{Remark}}

\newtheorem{Definition}[T]{Definition}
\newcommand{\cde}{\begin{Definition}}
\newcommand{\fde}{\end{Definition}}

\newcommand{\dem}{\begin{proof}}
\newcommand{\cqd}{\end{proof}}

\newtheorem{proposition}[theorem]{Proposition}
\newtheorem{lemma}[theorem]{Lemma}

\theoremstyle{definition}
\newtheorem{remark}[theorem]{Remark}

\newtheorem{definition}[theorem]{Definition}

\newcommand{\pf} {{\bf Proof: }}
\newcommand{\field}[1]{\mathbb{#1}}

\newcommand{\ov} {\overline}

\newcommand{\re}{\field{R}}

\newcommand{\integer}{\field{Z}}
\renewcommand{\natural}{\field{N}}

\newcommand{\al} {\alpha}       
\newcommand{\be} {\beta}

\newcommand{\la} {\lambda}      \newcommand{\La}{\Lambda}

\newcommand{\vsi}{\varsigma}

\newcommand{\SO}{{\mathbb O}}

\renewcommand{\SS}{{\mathbb S}}

\newcommand{\diam}{\operatorname{{diam}}}

\begin{document}

\author{Armando Castro}
\address{Departamento de Matem\'atica, Universidade Federal da Bahia\\
Av. Ademar de Barros s/n, 40170-110 Salvador, Brazil.}
\email{armandomat@pesquisador.cnpq.br}

\author{Krerley Oliveira}
\address{Departamento de Matem\'atica, Universidade Federal de Alagoas\\
Macei\'o, Brazil.}
\email{krerley@impa.br}

\author{Vilton  Pinheiro}
\address{Departamento de Matem\'atica, Universidade Federal da Bahia\\
Av. Ademar de Barros s/n, 40170-110 Salvador, Brazil.}
\email{viltonj@ufba.br}

\date{\today}

\thanks{Work carried out at the  Federal University of
Bahia. Partially supported by PADCT/CNPq and UFBA}

\title[Non Uniformly hyperbolic periodic points]{Shadowing by non
uniformly hyperbolic periodic points and uniform hyperbolicity }

\maketitle

\begin{abstract}
We prove that, under a mild condition on the hyperbolicity of its
periodic points, a map $g$ which is topologically conjugated to a
hyperbolic map (respectively, an expanding map) is also a hyperbolic
map (respectively, an expanding map). In particular, this result
gives a partial positive answer for a question done by A. Katok, in
a related context.
\end{abstract}

\section{Introduction}

Since Smale proposed the notion of \emph{uniformly hyperbolic
dynamical system}, the theory and results obtained by dynamicists
around the world have described many of its features, from the
structural and  measure-theoretical points of view.

Nevertheless,  the study of conditions for a non uniformly expanding
map be expanding is not well understood, regarding the few results
concerning the subject. One of these results, is the remarkable
theorem of Ma\~n\'e \cite{Ma}, valid for invariant sets without
critical points for interval maps. Outside this setting, not much is
known  and it is by itself a interesting point of research. In
particular, the study of non uniform expanding rates and conditions
over a  given set of points and its relations with uniform expanding
behavior appears in several recent papers (\cite{AAS}, \cite{YC} and
\cite{YLC}). Let us briefly describe some of this results:

We say that a {\em local diffeomorphism} $f$ is non uniformly
expanding ({\em NUE}) on a set $X$ , if there exists $\eta< 0$ such
that
$$
\liminf\limits_{n \rightarrow \infty} \frac 1 n\log
\sum_{j=0}^{n-1}\|[Df(f^j(x))]^{-1}\| \le \eta<0\mbox{ for all }x\in
X.$$ In (\cite{AAS}), the authors proved that any local
diffeomorphism in a compact manifold admitting non uniform expansion
at a set of total probability, i.e., with full measure for any
invariant measure, is in fact an expanding map. Similar results
holds  for diffeomorphisms.

By Oseledets(\cite{Os68}), one knows that if $\mu$ is an invariant
measure for a $C^1$ map $f$, then the number
$$
\lambda(x,v)= \lim\limits_{n\rightarrow \infty} \frac{1}{n} \log
\|Df^n(x)v\|,
$$
is defined in a set of total probability and it is called
\emph{Lyapunov exponent}  at $x$ in the direction $v$. In \cite{YC},
the author prove that if $f$ is a local diffeomorphism such all
Lyapunov exponents are positive then it is, in fact, a expanding map
and also obtained  the results for diffeomorphisms admitting
continuous splitting (see Theorems~\ref{t.cao1} and \ref{t.cao2}).

Here, we consider the universe of systems  without critical points
which are topologically conjugated to expanding maps. In such
context, a necessary and sufficient condition for a system to be
expanding is just that it is non uniformly expanding on the set of
periodic points. Therefore, as a main result, we prove that a local
diffeomorphism topologically conjugated to an expanding map is
itself an expanding map if, and only if, it is non uniformly
expanding on the set of periodic points. We also obtain a similar
result for dynamics with non uniformly hyperbolic (NUH) periodic
points conjugated to an uniformly hyperbolic map.

\cmt \label{theo1}  Let $g:M \to M$ be a $C^2$-class local
diffeomorphism on a compact manifold $M$. Suppose that $g$ is
topologically conjugated to an expanding $C^1$ map $f$. If $g$ is
non uniformly expanding on the set $Per(g)$ of periodic points, then
$g$ is an expanding map. \fmt

\begin{figure}
\label{fig1}
\begin{center}
\includegraphics{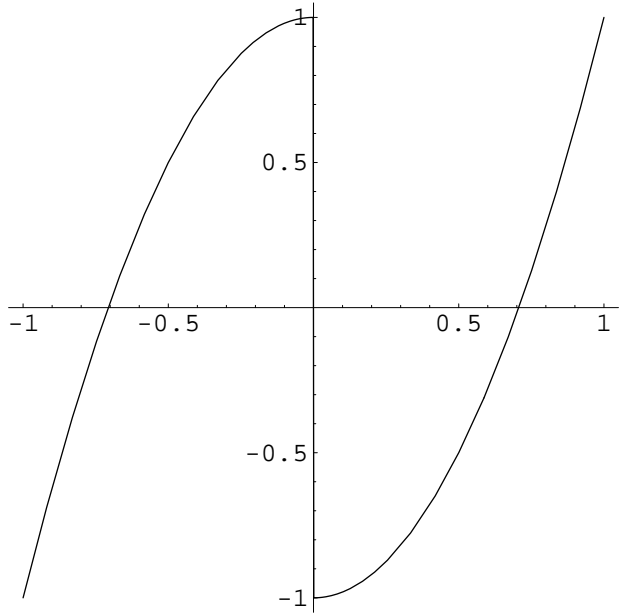}
\end{center}
\caption{Lift to the interval of a map in $S^1$ satisfying NUE
condition on the periodic points and topologically conjugated to $z
\mapsto z^2$}
\end{figure}

\begin{remark}\label{r.pontoscriticos} We observe that the condition NUE
on the periodic points  is not enough to assure that the map $g$ is
expanding, even if we assume that $g$ is topologically conjugate to
an expanding map. It is a standard matter that the map $z
\rightarrow z^2$, defined on the circle is topologically conjugated
to a map with criticalities  satisfying NUE condition on the
periodic points. See the figure \ref{fig1}.
\end{remark}

In the Theorem~\ref{theo1}, due to the fact that we are dealing with
maps that are local diffeomorphisms we avoid examples as in
Remark~\ref{r.pontoscriticos}.

For diffeomorphisms, the existence of a continuous splitting of $M$
plays a similar role.  In \cite{YLC}, the authors exhibit an example
of a non-hyperbolic horseshoe such that the splitting is continuous
over the periodic points and all Lyapunov exponents are positive and
bounded from zero. In particular, some condition of continuity of
the splitting in the closure of the periodic points is necessary. In
order to state our results in the invertible case, we need the
following definition:

\begin{definition}{(Non uniformly hyperbolic set).} \label{defNUH}
\label{defnuh} Let $g: M \to M$ be a diffeomorphism on a compact
manifold $M$. We say that an invariant set $S \subset M$ is a {\em
non uniformly hyperbolic set} or, simply, {\em NUH}, iff
\begin{enumerate}
\item
There is an $Dg-$invariant splitting $T_S M= E^{cs} \oplus E^{cu}$;

\item There exists $\eta< 0$ and an
adapted Riemannian metric for which any point $p \in S$  satisfies
$$
\liminf_{n \to +\infty}\frac{1}{n}\sum_{j= 0}^{n-1} \log
\|Dg(g^j(p))|_{E^{cs}(g^j(p))}\| \leq \eta
$$
and
$$
\liminf_{n \to +\infty}\frac{1}{n}\sum_{j= 0}^{n- 1} \log
\|[Dg(g^{j}(p))|_{E^{cu}(g^{j}(p))}]^{-1}\| \leq \eta
$$

\end{enumerate}

\end{definition}

We also recall here the notion of hyperbolic set:
\begin{definition}
Let $\La$ be an invariant set for a $C^1$ diffeomorphism $f$ of a
manifold $M$. We say that $\La$ is a hyperbolic set if there is a
continuous splitting  $T_\La M= E^s \oplus E^u$ which is
$Tf$-invariant ($Tf(E^s)= E^s, Tf(E^u)= E^u$) and for which there
are constants $c> 0$, $0< \la< 1$, such that
$$
\|Tf^n|_{E^s}\| < c \la^n, \quad  \|Tf^{-n}|_{E^u}\| < c \la^n,
\forall n \in \natural.
$$
\end{definition}

For the diffeomorphism case, we have two slightly different results:

\cmt \label{theo2} Let $g:M \to M$ be a $C^2$ diffeomorphism on a
compact manifold $M$, and let $\La \subset M$ be a compact invariant
set. Suppose that $g|_{\La}$ is topologically conjugated to a $C^1$
diffeomorphism $f$ restricted to a set $\hat \La$, hyperbolic for
$f$. If the set $Per(g)$ of periodic points of $g$ is non uniformly
hyperbolic (NUH), and $T_{Per(g)}M= E^{cs} \oplus E^{cu}$ is a
dominated splitting, then $\La$ is a hyperbolic set for $g$. \fmt

\cmt \label{theo3}  Let $g:M \to M$ be a $C^2$ diffeomorphism on a
compact manifold $M$, and let $\La \subset M$ be a compact invariant
set. Suppose that $g|_{\La}$ is topologically conjugated to a $C^1$
diffeomorphism $f$ restricted to a set $\hat \La$, hyperbolic for
$f$. If the set $Per(g)$ of periodic points of $g$ is non uniformly
hyperbolic (NUH), and $T_{Per(g)}M= E^{cs} \oplus E^{cu}$ has a
continuous extension to a splitting on $T_{\ov{Per(g)}}M$, then
$\La$ is a hyperbolic set for $g$. \fmt

In fact, theorem \ref{theo2} is a consequence of this last theorem
\ref{theo3}. Nevertheless, the hypotheses in \ref{theo2} are easier
to verify.

\begin{remark}
With the same proof, all results in this paper are valid if the
derivative $Dg$ is just H\"older continuous.
\end{remark}

\begin{definition}{(Shadowing by periodic points).} Let $f:M \to M$
be a map and $\hat\La \subset M$ be a compact $g$-invariant set. We
say that $(f,\hat\La)$ has the {\em shadowing by periodic points
property} if given $\epsilon> 0$, exists $\al> 0$ such that for any
orbit segment $\{x, \dots, f^n(x)\} \subset \hat\La$ with $d(f^n(x),
x)< \al$ there exists a periodic point $p \in M$ with period $n$
such that $d(f^j(p), f^j(x))< \epsilon, forall 0 \leq j \leq n$. In
this case, we say that the orbit of $p$ {\em $\epsilon-$shadows the
orbit segment} $\{x, \dots, f^n(x)\}$.
\end{definition}

If $\hat\La \subset M$ is a hyperbolic, compact invariant set for a
diffeomorphism $f$, then the classical theory of hiperbolic systems
implies that $(f,\hat \La)$ has the {\em shadowing by periodic
points property} (see proposition 8.5 in \cite{Sh87}). The same is
also valid for any system which is topologically conjugated to $f$.
{\em Shadowing by periodic points} is the key ingredient in the
proofs of the theorems \ref{theo1}, \ref{theo2}, and \ref{theo3} we
stated above in this introduction. Therefore, as a consequence of
their proofs, we also obtain the following (more general) results:

\cmt \label{theo4}  Let $g:M \to M$ be a $C^2$ local diffeomorphism
on a compact manifold $M$. Suppose there exists an invariant compact
set $\La \subset M$ such that $(g, \La)$ has the {\em shadowing by
periodic points property} If $g$ is non uniformly expanding on the
set $Per(g)$ of periodic points, then $g$ is an expanding map on
$\La$.
\fmt

\cmt \label{theo5} Let $g:M \to M$ be a $C^2$ diffeomorphism on a
compact manifold $M$,  and let $\La \subset M$ be a compact
$g$-invariant set. Suppose that $(g, \La)$ has the {\em shadowing by
periodic points property}. If the set $Per(g)$ of periodic points of
$g$ is non uniformly hyperbolic (NUH), and $T_{Per(g)}M= E^{cs}
\oplus E^{cu}$ is a dominated splitting, then $\La$ is a hyperbolic
set for $g$. \fmt

\cmt \label{theo6} Let $g:M \to M$ be a $C^2$ diffeomorphism on a
compact manifold $M$, and let $\La \subset M$ be a compact
$g$-invariant set. Suppose that $(g, \La)$ has the {\em shadowing by
periodic points property}. If the set $Per(g)$ of periodic points of
$g$ is non uniformly hyperbolic (NUH), and $T_{Per(g)}M= E^{cs}
\oplus E^{cu}$ has a continuous extension to a splitting on
$T_{\ov{Per(g)}}M$, then $\La$ is a hyperbolic set for $g$. \fmt

\section{The endomorphism case: non uniformly expanding periodic set}

During this section, $g:M \to M$ will always be a $C^2-$local
diffeomorphism which is topologically conjugated to a $C^1$
expanding endomorphism.

We recall the definition of NUE:

\begin{definition}{(Non uniformly expanding set).} \label{defNUE}
\label{defnue} Let $g: M \to M$ be a map on a compact manifold $M$.
We say that an invariant set $S \subset M$ is a { \em non uniformly
expanding set} or, simply, {\em NUE}, iff:

There exists $\eta< 0$ and an adapted Riemannian metric for which
any point $p \in S$  satisfies
$$
\liminf_{n \to +\infty}\frac{1}{n}\sum_{j= 0}^{n-1} \log
\|[Dg(g^j(p))]^{-1}\| \leq \eta
$$

\end{definition}

\begin{remark}\label{reNUE}
In this paper, we will always focus our attention on the set of
periodic points $Per(g)$ of $g$. Given a point $p \in Per(g)$, let
us set $t= t(p):= period(p)$. In such case, the following
equivalence is immediate:
 $S:= Per(g)$ is NUE iff there exists $\vsi< 1$ such that for each periodic point $p$, $\prod_{j=
0}^{t(p)-1} \|[Dg(g^j(p))]^{-1}\|< \vsi^{t(p)}$.
\end{remark}

For the sequence, we give a simplified definition for the case of local
diffeomorphisms of the notion introduced by \cite{ABV00}:

\begin{definition}{(Hyperbolic time for local diffeomorphisms)}
Let  $z \in M$ be a regular point. We say that $k \in \natural$ is a
{\em $\vsi$-hyperbolic time} for $z$ if for  $i= 1, \dots, k$, holds
\begin{equation}
\prod_{j= 1}^{i}\|[Dg|(g^{p-j}(z))]^{-1}\| \leq \vsi^i.
\end{equation}

\end{definition}

\begin{lemma}\label{leNUE}
Suppose that $g$ is topologically conjugated to an expanding map
$f$. Let $x$ be a recurrent, regular point of $g$. If $Per(g)$ is
NUE, then all Lyapunov exponents of $x$ are positive.
\end{lemma}
\pf Let $\delta> 0$ such that, given any ball $B(z, \delta)$ the
corresponding inverse branches of $g$ are well defined
diffeomorphisms. Let $\vsi= e^\eta$, $\eta$ as in definition
\ref{defnue}, $\vsi< \vsi' < 1$ fixed,  and let $\epsilon> 0$ such
that $(\sqrt{\vsi'})^{-1}- \epsilon > 1$. Since $x$ is a regular
point, there is $n_0 \in \natural$ such that
$$
(\vsi_j- \epsilon)^n\cdot \|v_j\| <
\|Dg^n(x) \cdot v_j\| < (\vsi_j+ \epsilon)^n \cdot \|v_j\|
\forall v_j \in E_j, \forall n \geq n_0.
$$
where $E_j$ are the Lyapunov eigenspaces and $\log(\vsi_j)$ are
their respective Lyapunov exponents.

Now, by Pliss Lemma~\cite{Pli}, there exists $n_1 > n_0$ such that
any point $y$ for which we have $\prod_{j=
0}^{n-1}\|[Dg(g^j(y))]^{-1}\|^{-1} \geq \vsi^{-n}$, for some $n \geq
n_1$, then $y$ has, at least, $n_0$ $\vsi'-$hyperbolic times less
than $n_1$.

We fix $0< \delta' \leq \delta$ such that
$$
\|[Dg^{-1}(y)\| \leq \frac{1}{\sqrt{\vsi'}} \|Dg^{-1}(z)\|, \forall
z, y; d(z, y)< \delta',
$$
where $g^{-1}$ is an inverse branch for $g$.

We set $0< \delta''< \delta'$ such that if $g^{-n}$ is an arbitrary
composition of $n$ inverse branches for $g$, then $\diam( g^{-n}
(B(z, \delta''))< \delta'$, $\forall z \in M$, $\forall n \in
\natural$. This occurs because it is valid for the hyperbolic system
$f$ to which $g$ is conjugated.

As $x$ is a recurrent point, we set $n_2 \geq n_1$ a return time
such that a neighborhood $V_x \subset B(x, \delta'')$ of $x$ is
taken by $g^{n_2}$ onto $B(x, \delta'')$.

Therefore, writing $G:= (g^{n_2}|_{V_x})^{-1}$,
$G: B(x, \delta'') \to V_x \subset B(x, \delta'')$ has a
fixed point $p \in V_x$, which is a
periodic point of period $n_2$ for $g$. By hypothesis,
$p$ is a hyperbolic periodic point for which we have
$$
\prod_{j= 0}^{n_2-1}\|[Dg(g^j(p))]^{-1}\|^{-1}\geq \|\vsi^{-n_2}\|
\Rightarrow \|DG(p)\| \leq \|\vsi^{n_2}\|.
$$

By our choice of $n_1$ and the equation above, there
exists a $\vsi'-$hyperbolic time $n_0< n'< n_2$ for
$p$.

Due to lemma 2.7 in \cite{ABV00} (see also prop. 2.23 in
\cite{C2}), $n'$ is also a $\sqrt{\vsi'}-$hyperbolic
time for $x$. In particular, this implies that
$$
\|Dg^{n'}(x) \cdot v\| \geq \sqrt{\vsi'}^{-n'} \|v\|,
\forall v \in T_p M.
$$
Therefore, $\vsi_j\geq \sqrt{\vsi'}^{-1}- \epsilon> 1$,
$\forall j$. This means that all Lyapunov exponents of
$x$ are greater than $1$.

\qed

We note that the set of Oseledet's regular, recurrent points is a
total probability set, due to Oseledet's theorem and Poincar\'e's
Recurrence theorem. This means that such set has measure equal to
$1$ for any $g-$invariant probability measure. So, for any
$g$-invariant measure, lemma \ref{leNUE} implies that all Lyapunov
exponents are positive.
Therefore, our theorem \ref{theo1} is obtained applying lemma
\ref{leNUE} to the following result:

\begin{theorem}{\cite{YC}}\label{t.cao1}
Let $f: M \to M$ be a $C^1$ local diffeomorphism on a compact
Riemannian manifold. If the Lyapunov exponents of every $f$
invariant probability measure are positive, then $g$ is uniformly
expanding.
\end{theorem}

\section{The diffeomorphism case: non uniformly hyperbolic periodic
set}

Now, we treat the case when $f$ is a diffeomorphism. Along this
section, we suppose that the periodic set $Per(g)$ is NUH. (see
definition \ref{defNUH}, on page \pageref{defNUH}).

The following remark is the analogous of remark \ref{reNUE} for the
diffeomorphism case:

\begin{remark}
We note that the set of period points $Per(g)$ is NUH iff there
exists $\vsi> 1$ such that for each periodic point $p$ with period
$t(p)$, then $\prod_{j= 0}^{t(p)-1}
\|[Dg|_{E^u}(g^j(p))]^{-1}\|^{-1}> \vsi^{t(p)}$ and $\prod_{j=
0}^{t(p)-1} \|Dg|_{E^s}(g^j(p))\|< \vsi^{-t(p)}$.
\end{remark}

Before we state and prove the next lemma let us
introduce some notation. Given a periodic
point $p \in M$, we
denote the cone over $E^s(p)$ of width $0< a< 1$
by
$$
C^s_a(p):= \{v_s+  v_u \in E^s(p) \oplus  E^u(p),
\text{ such that } a \| v_s\|< \|v_u\|.
$$
Analogously, we define a cone over $E^u(p)$
of width $a$.

Now we adapt the definition of hyperbolic times to the context of
diffeomorphisms (see also \cite{C1}).

\begin{definition}{(Hyperbolic time for stable directions)}
 Let $0< \la < 1$ and $z \in M$ be a regular point.
 Suppose that $E$ is an invariant subbundle of $T_{\SS(z)}M$,
 where $\SS(z)$ is some  orbit segment of $z$.
 We say that $k \in \natural$ is
a {\em $\la-$hyperbolic time} for $z$ if for $g^{-k}(z)= y$ and $i=
1 \dots k$, holds
\begin{equation}
\prod_{j= 0}^{i-1}\|Dg|_E(g^j(y))\| \leq \la^i.
\end{equation}
\end{definition}

An analogous definition can be done for unstable directions just
exchanging $g$ by $g^{-1}$ in the definition above.

\begin{lemma}\label{leNUH}

Let $g:M \to M$ be a $C^2$ diffeomorphism and $\Lambda \subset M$ be
some compact $g$-invariant set. Suppose that $g|\Lambda$ is
topologically conjugated to a  $f|_{\hat \Lambda}$, where $\hat
\Lambda$ is a hyperbolic set for $f$. Let $x$ be a recurrent,
regular point of $g$. Suppose that $Per(g)$ is NUH, and that the
splitting $T_{Per(g)}= E^{cs} \oplus E^{cu}$ have a continuous
extension to $T_{\ov{Per(g)}}M= E^1 \oplus E^2$. Then all Lyapunov
exponents of $x$ are nonzero.
\end{lemma}
\pf Let $\vsi= e^\eta$, $\eta< 0$ as in definition \ref{defnuh},
$\vsi< \vsi' < 1$ fixed,  and let $\epsilon> 0$ such that
$(\sqrt{\vsi'})^{-1}- \epsilon > 1$. Since $x$ is a regular point,
there is $n_0 \in \natural$ such that
$$
(\vsi_j- \epsilon)^n\cdot \|v_j\| < \|Dg^n(x) \cdot v_j\| < (\vsi_j+
\epsilon)^n \cdot \|v_j\| \forall v_j \in E_j, \forall n \geq n_0,
$$
and
$$
(\vsi_j- \epsilon)^{-n}\cdot \|v_j\| >\|Dg^{-n}(x) \cdot v_j\| >
(\vsi_j+ \epsilon)^{-n} \cdot \|v_j\| \forall v_j \in E_j, \forall n
\geq n_0.
$$
 where $E_j$ are the Lyapunov eigenspaces and $\log(\vsi_j)$ are
their respective Lyapunov exponents. We denote by $E^cs(x)$
(respectivelly,  $E^{cu}(x)$)  the space spanned by the Lyapunov
eigenspaces with negative (respectivelly, positive) Lyapunov
exponents. $E^0(x)$ will denote the Lyapunov eigenspace
corresponding to an eventual zero Lyapunov exponent.

Let us prove that the dimension of the space $E^{cs}(x)$
corresponding to the negative Lyapunov exponents of $x$ is equal or
greater than the dimension of the stable space of any periodic
point. An analogous result will obviously hold for $E^{cu}(x)$.
Therefore, we conclude that $T_x M= E^{cs}(x) \oplus E^{cu}(x)$ and
that all Lyapunov exponents of $x$ are nonzero.

By taking charts, due to the uniform continuity of $Dg$, we can fix
$0< a < 1$ and $0< \delta'$ such that if $z$ is periodic,
$$
\|Dg(y) \cdot v\| \leq \frac{1}{\sqrt{\vsi'}}
\|Dg|_{E^{cs}(z)}(z)\|\|v\|, \forall y \in B(z,\delta'), \forall v
\in C^s_a(z).
$$
Due to the continuity of $E^1$, we can assume $a$ small enough such
that each cone $C^s_a(z)$ contains $E^1(y)$ or all point $y \in B(z,
\delta') \cap \ov{Per(g)}$.

Now, by Pliss Lemma in \cite{Pli}, there exists $n_1 > n_0$ such
that any point $z \in {Per(g)}$ for which we have
$\prod_{j=0}^{n-1}\|Dg|_{E^1}(g^{-n+ j}(z)\| \leq \vsi^{n}$,  for
some $n \geq n_1$, then $z$ has, at least, $n_0$ $\vsi'-$hyperbolic
times less than $n_1$.

As $x$ is a recurrent point (also for $g^{-1}$), we set $n_2 \geq
n_1$ a return time for $g^{-1}$ such that there exists one periodic
point $p$ with period $n_2$ that $\delta'/3-$shadows the orbit
segment $\{x, g^{-1}(x), \dots, g^{-n_2}(x)\}$. Such periodic point
exists because $g|_\Lambda$ is conjugated to a diffeomorphism
$f|_{hat \Lambda}$ which is shadowed by periodic points (see
proposition 8.5 in \cite{Sh87}).

By hypothesis, $p$ is a hyperbolic periodic point for which we have
$$
\|\prod_{j= 0}^{n_2-1} Dg|_{E^{cs}}(g^j(p))\|\leq \|\vsi^{n_2}\|
$$

By our choice of $n_1$ and the equation above, there exists a
$\vsi'-$hyperbolic time $n_0< n'< n_2$ for $p$.

Due to  prop. 2.23 in \cite{C2}, $n'$ is also a
$\sqrt{\vsi'}-$hyperbolic time for $x$. More precisely, this means
that
$$
\prod_{j= 0}^{n'-1}\| Dg|_{E^1(g^{-n'+ j}(x))}(g^j(g^{-n'}(x)))\|
\leq \sqrt{\vsi'}^j,
$$
since the space  $E^1(g^{-n'+ j}(x)) \subset C^s_a(g^{-n'+j}(p))$.

In particular, this implies that

\begin{equation}
\|Dg^{-n'}(x)\cdot v\| \geq \sqrt{\vsi'}^{-n'} \|v\|, \forall v \in
E^1(x) \label{eqarremat}
\end{equation}

This implies that the dimension of the negative Lyapunov exponents
space $E^{cs}(x)$ is at least the dimension of $E^1(x)$ which equals
the dimension of $E^{cs}(p)$. In fact, if we had $\dim(E^{1}(x))>
\dim(E^{cs}(x))$, then since $T_xM= E^{cs}(x) \oplus E^0(x) \oplus
E^{cu}(x)$, the intersection $E^1(x) \cap  (E^0(x) \oplus
E^{cu}(x))$ would be nontrivial. That is an absurd, because no
vector in $(E^0(x) \oplus E^{cu}(x)) \setminus \{ 0\}$ satisfies
equation \ref{eqarremat}.

Applying the same arguments above to $E^{cu}(x)$ we conclude that
the number of positive Lyapunov at $x$ is at least the dimension of
$E^cu(p)$ and this concludes the lemma.

\qed

For the next results we recall here the definition of {\em dominated
splitting}:

\begin{definition}{\em(Dominated splitting).} \label{defdom}
Let $f:M \to M$ be a diffeomorphism on a compact manifold $M$ and
let $X \subset M$ be an invariant subset. We say that a splitting
 $T_X M= E \oplus \hat E$ is a {\em dominated splitting} iff:
\begin{enumerate}
\item The splitting is invariant by $Df$, which means that
$Df(E(x))= E(f(x))$ and $Df(\hat E(x))= \hat E(f(x))$.
\item There exist $0< \la< 1$ and some $l \in \natural$ such that
for all $x \in X$
$$
\sup_{v \in E, \|v\|= 1} \{\|df^l(x) v\|\} \cdot (\inf_{v \in \hat
E, \|v\|= 1} \{\|df^l(x) v\|\})^{-1} \leq \la.
$$

\end{enumerate}

\end{definition}

A priori, we do not require dominated splitting to be continuous.
However, they always are:

\begin{lemma}\label{ledom}
Let $f:M \to M$ be a diffeomorphism on a compact manifold $M$. Let
$X\subset M$ be some $f-$invariant set. Suppose there exists some
invariant dominated splitting $T_X M= E \oplus \hat E$. Then, such
splitting is continuous in $T_X M$, and unique since we fix the
dimensions of $E, \hat E$. Moreover, it extends uniquely and
continuously to a splitting of $T_{\ov X}M$.
\end{lemma}
\pf
By replacing $f$ by an iterate, there is no loss of generality
in supposing that $l \in \natural$ in definition \ref{defdom} equals
to 1.
We start by constructing an invariant dominated splitting on
$T_{\ov X} M$ extending the one we have on $T_X M$. Let $\SO(x)$ be
an orbit contained in $\ov X$. Our construction will be dependent of
some choices. We choose one representative of $\SO(x)$, for example
$x$. Let us also choose some $(x_n)$, $x_n \in X$, $x_n \to x \in
M$, as $n \to \infty$. Let $v^1_n, \dots, v^s_n \in E(x_n)$, $\hat
v^{s+1}_n, \dots, \hat v^m_n \in \hat E(x_n)$ be orthonormal bases
of $E(x_n)$, $\hat E(x_n)$, respectively. The domination property is
equivalent to
$$
\|Df(x_n) \sum_{j= 1}^s \al_j v^j_n\|\cdot \|Df(x_n) \hat \sum_{i=
s+1}^m \be_i \hat v^i_n\|^{-1} \leq \la < 1,
$$
for any convex combination $\sum_{j= 1}^s \al_j v^j_n$, $\sum_{i=
s+1}^m \be_i v^i_n$.  Replacing by some convergent subsequence, if
necessary, we can suppose that $(v^1, \dots,$ $v^s)$, $v^1, \dots, v^s
\in T_x M$ (resp. $(\hat v^{s+1}, \dots, \hat v^m)$ ) is the limit
of the sequence $(v^1_n, \dots, v^s_n)$ (resp. of the sequence
$(\hat v^{s+1}_n, \dots, \hat v^m_n)$). Since the domination
property is a closed condition,
$$
\|Df(x) \sum_{j= 1}^s \al_j v^j\|\cdot \|Df(x) \hat \sum_{i= s+1}^m
\be_i \hat v^i\|^{-1} \leq \la < 1,
$$
holds.

Now, we write $G$ for the Gram-Schmidt operator (which takes a
linearly independent set of vectors on an orthonormal set of vectors
spanning the same vector space). Given any iterate $y= f^k(x), k \in
\integer$, then $f^k(x_n) \to y$ and
$$
G \circ (Df^k(x_n)v^1_n, \dots, Df^k(x_n)v^s_n) \to G \circ
(Df^k(x)v^1, \dots, Df^k(x)v^s),
$$
$$
G \circ (Df^k(x_n)\hat v^{s+1}_n, \dots, Df^k(x_n)\hat v^m_n) \to G
\circ (Df^k(x)\hat v^{s+1}, \dots, Df^k(x)\hat v^m),
$$
as $n \to \infty$.  Writing $(w^1_k, \dots, w^s_k):= G \circ
(Df^k(x)v^1, \dots, Df^k(x)v^s)$ and  $ (\hat w^{s+1}_k, \dots, \hat
w^m_k):= G \circ (Df^k(x)\hat v^{s+1}, \dots, Df^k(x)\hat v^m)$, $k
\in \integer$, the same calculations above show that
$$
T_y M= span\{w^1_k, \dots, w^s_k\} \oplus span\{\hat w^{s+1}_k,
\dots, \hat w^m_k\}=: E(y) \oplus \hat E(y)
$$
is a dominated splitting.

Moreover, it is clear that
$$
Df(f^k(x)) (span\{w^1_k, \dots, w^s_k\})= span\{w^1_{k+1}, \dots,
w^s_{k+1}\}
$$
and
$$
Df(f^k(x)) (span\{\hat w^{s+1}_k, \dots, \hat w^m_k\})=
span\{w^{s+1}_{k+1}, \dots, \hat w^m_{k+1}\},
$$
which implies that it is an invariant splitting.

Note  that since the dominated splitting condition is a closed
condition, if we prove that there  exists  a unique dominated
splitting with the same dimensions of the splitting we constructed,
it will be automatically continuous. This is because, given $x_n \to
x \in X$, any convergent sequences of orthonormal bases of
$E(x_n)$,$\hat E(x_n)$ will converge to  orthonormal bases of
dominated spaces in $T_x M$ which due the uniqueness will
necessarily be $E(x)$, $\hat E(x)$.

The argument to prove uniqueness is the following. Suppose that we
have two invariant dominated splittings $T_{\ov X} M= E \oplus \hat
E$, $T_{\ov X} M= E' \oplus \hat E'$. Fix an arbitrary $x \in \ov
X$.

By changing $f$ for some positive iterate $f^l$ , there is no loss
of generality in supposing that domination condition is valid for
$l= 1$ on both splittings. The domination condition yields:
$$
\|df|_{E(x)}\|(\inf_{v \in \hat E, \|v\|=1} \{\|Df(x)v\|\})^{-1}
\leq \la
$$
and
$$\|df|_{E(x)}\|(\inf_{v \in \hat E',
\|v\|=1} \{\|Df(x)v\|\})^{-1} \leq \la.
$$
Let us show that $E(x)= E'(x)$. Note that if $E(x) \subset E'(x)$
(or vice-versa), as the spaces have the same dimension, they should
be the same. So, let us suppose by contradiction that there exist $v
\in E(x) \setminus E'(x)$ and $v' \in E'(x) \setminus E(x)$. We then
write $v= v_{E'} + v_{\hat E'}$, with $v_{E'} \in E'$, $v_{\hat E'}
\in \hat E'$ and $v_{\hat E'} \neq 0$. This last inequality,
together with the invariance of the splittings implies that
$$
Df^n(x) \cdot v= \al_n  v^n_{E'}+ \be_n v^n_{\hat E'},
$$
where $v^n_{E'}$ and $v^n_{\hat E'}$ are unitary vectors
respectively in $E'(f^n(x))$ and $\hat E'(f^n(x))$ and $\al_n/\be_n
\lessapprox \la^n \to 0$. In particular, $Df^n(x) \cdot v \in
E(f^n(x))$ belongs in an arbitrarily small width cone over $\hat
E'(f^n(x))$ (which dominates $E'(f^n(x))$), as we take $n$
sufficiently big. This implies that, for all $n \in \natural$
sufficiently big, there exists $v_n= Df^n(x) \cdot v/\|Df^n(x) \cdot
v\| \in E(f^n(x))$ such that
$$
\|df|_{E'( y_n)}\| \cdot \|Df( y_n) v_n\|^{-1} < \tilde \la < 1,
$$
where $y_n= f^n(x)$. Now, we repeat the same argument above for $v'
\in E'(x) \setminus E(x)$ and again for all $n$ sufficiently big, we
obtain unitary vectors $v'_n \in E'(y_n)$ such that
$$
\|df|_{E( y_n)}\| \cdot \|Df( y_n) v'_n\|^{ -1} < \tilde \la < 1.
$$
Therefore, we have
$$
\|df|_{E'( y_n)}\| \leq \tilde \la \cdot \|Df( y_n) v_n\| \leq
\tilde \la \cdot \|df|_{E( y_n)}\|
$$
and
$$
\|df|_{E( y_n)}\| \leq \tilde \la \cdot \|Df( y_n) v'_n\| \leq
\tilde \la \cdot \|df|_{E( y_n)}\|,
$$
which is a contradiction.

 \qed

\begin{lemma}\label{leNUH2}
Suppose that $g$ is topologically conjugated to a hyperbolic map
$f$. Let $x$ be a recurrent, regular point of $g$. Suppose that
$Per(g)$ is NUH, and that the splitting $T_{Per(g)}M= E^{cs} \oplus
E^{cu}$ is a dominated splitting. Then all Lyapunov exponents of $x$
are nonzero.
\end{lemma}
\pf The proof is a direct consequence of lemmas \ref{leNUH} and
\ref{ledom}. By lemma \ref{ledom}, the invariant dominated splitting
over $T_{Per(g)}M$ extends to a unique continuous invariant
(dominated) splitting over $\ov T_{Per(g)}M$. So, we become under
the hypotheses of lemma \ref{leNUH}, which allows to conclude that
all Lyapunov exponents of any recurrent point $x \in M$ are nonzero.

\qed

By the same arguments as in the expanding case (see paragraph below
the proof of lemma \ref{leNUE}), our theorem \ref{theo2} is obtained
applying lemma \ref{leNUH} to the following result:

\begin{theorem}{\cite{YC}}\label{t.cao2}
Let $f: M \to M$ be a $C^1$ diffeomorphism on a compact Riemannian
manifold, with a positively invariant set $\La$ for which the
tangent bundle has a continuous splitting $T_\La M = E^{cs} \oplus
E^{cu}$. If $f$ has positive Lyapunov exponents in the $E^{cu}$
direction and negative Lyapunov exponents in the $E^{cs}$ direction
on a set of total probability, then $f$ is uniformly hyperbolic.
\end{theorem}

\section{On a conjecture of A. Katok}

A. Katok has conjectured that a $C^{1+}$ system which is H\"older
conjugated to an expanding map  (respectively, an Anosov
diffeomorphism) is also expanding (respectively, is also an Anosov
diffeomorphism).

Note that, under the hypotheses of such conjecture, the periodic points
of the $g: M \to M$ are hyperbolic,
with uniform
bounds for the eigenvalues of iterate of $Dg$ in the period of such points.
This is proven below.

First, we consider the expanding case. Let $p$ a periodic
point of period $t$ of $g$. Then, $h(p)$ is a periodic
point of period $t$ of $f$.
Let us call $f^{-1}$ the inverse
branch of $f$, defined on a neighborhood of the orbit
of $h(p)$, for which $h(p)= \hat p$ is a periodic point of
period $t$.
Analogously, let us call $g^{-1}$ be the inverse
branch of $g$ for which $p$ is a periodic point of
period $t$. Since $f$ is an expanding map, there are
$0< \hat \la < 1$ and $\hat \delta> 0$ such that
$$
d(f^{-j}(\hat x), f^{-j}(\hat y))\leq{\hat\la}^j d(\hat x,\hat y),
\forall j \in \natural,
\forall \hat x, \hat y \in B(\hat p, \hat \delta).
$$
As an immediate consequence of the $C^\alpha$ conjugation
$h$ there exists $\delta> 0$ such that
$$
d(g^{-j}(x), g^{-j}(y))\leq({\hat\la}^\al)^j K^{1+\al} d(x, y)^{(\alpha^2)},
\forall j \in \natural,
\forall x, y \in B(p,  \delta)
$$
and
\begin{eqnarray}
d(g^{-j}(x), g^{-j}(y))\leq({\hat\la}^\al)^j K^{1+ \al} \delta^{\alpha^2},
\forall j \in \natural,
\forall x, y \in B(p,  \delta).
\label{eq1}
\end{eqnarray}

\begin{proposition}
Let  $B(x_0, r) \subset M$ and $G:\ov{B(x_0, r)} \to B(x_0, r)$ a
class $C^1$ local diffeomorphism such that $G(x_0)= x_0$ and for
some $0< \la < 1$ and $0< \be< 1$
$$
d(G^{n}(x), G^{n}(y)) \leq \la^n d(x, y)^\be, \forall
x, y \in B(x_0, r).
$$

Then all eigenvalues of $DG(x_0)$ are equal or less than $\la$.

\end{proposition}
\pf
Using charts, there is no loss of generality in supposing
that $M$ is an euclidean space and $x_0= 0$.
By contradiction, suppose there exists an invariant splitting
$\re^m= E^s+ E^c$, an adapted norm
$\|x\|= \|(x_s, x_c)\|= \max\{\|x_s\|, \|x_c\|\}$
and $\sigma> \la$ such that
$$
\|DG(0)\cdot x_s\| \leq  \la \cdot \|x_s\|, \forall x_s \in
E^s,
$$
$$
\|DG(0)\cdot x_c\| \geq  \sigma \cdot \|x_c\|, \forall x_c \in
E^c.
$$
Let $\epsilon> 0$ such that $\la+\epsilon< \sigma- \epsilon$
and take $\theta= \frac{\la+\epsilon}{\sigma- \epsilon}$.

Therefore, there is $\tilde r \leq r$ such that if we write
$$
G(x)= DG(0)\cdot x + \rho(x),
$$
then $\|\rho(x)\| < \epsilon \|x\|$,
$\forall x, \|x\|< \tilde r$.

We define a central cone
$$
V_c:= \{(x_s, x_c); \|x_s\| \leq \theta \|x_c\|\}
$$
By the hypothesis, there exists $\tilde r \leq \tilde r$
such that $G^n(B(0, \tilde r)) \subset B(0, \tilde r),$ $\forall n \in \natural$.
So, let us iterate $x \in B(0, \tilde r) \cap V_c$ (we
write $x^n= G^n(x)$). We obtain:
$$
\|x^1_c\| \geq \sigma \|x^0_c\|- \epsilon \|x^0\|
\geq (\sigma- \epsilon) \|x^0_c\|
$$
and
$$
\|x^1_s\| \leq \la \|x^0_s\|+ \epsilon \|x^0\| \leq (\la+ \epsilon) \|x^0_c\|.
$$
This implies that
$$
\|x^1_s\| \leq \frac{\la+\epsilon}{\sigma- \epsilon} \|x^1_c\|.
$$
In particular, if $x \in B(0, \tilde r) \cap V_c$ then
$G(x) \in V_c$.

% \|x^1_s\|/\|x^1_c\| \leq \theta
Therefore proceeding inductively, we obtain
$$
\|x^n\|= \|x^n_c\| \geq (\sigma- \epsilon)^n \|x^0_c\|= (\sigma- \epsilon)^n \|x^0\|.
$$
This contradicts the hipothesis, which implies that
$$
\|x^n\| \leq const \cdot \la^n, \forall n \in \natural.
$$

As $\epsilon>0 $ is arbitrary, we conclude that any eigenvalue of
$DG(0)$ is less than $\la$.

\qed

The proposition above implies our assertive that, if a map $g$ is
H\"older conjugated to an expanding map (respectively, Anosov) then
all periodic points have only nonzero Lyapunov exponents, and such
exponents are uniformly bounded away from zero. However, up to now
we do not know if, for example, the mild uniformity given by a
H\"older conjugation, plus the conjugation itself, imply that
$Per(g)$ is NUE.

Nevertheless, as a direct consequence of the last section, we obtain
that such conjecture is valid in the case that $Per(g)$ is NUE
(respectively, for Anosov, if $Per(g)$ is NUH with dominated
splitting).

\end{document}